# A stable algorithm for non-homogeneous waveguide equation based on DtN maps

Yin Wang, Jinyang Huang[1]

College of Science, Beijing University of Chemical Technology, CHINA

**Abstract:** A new stable computational method for non-homogeneous waveguide equation with a piecewise uniform structure along the main propagation direction is constructed, based on the modified Dirichlet-to-Neumann (DtN) map of each uniform segment. For segments with the same structure, only a DtN map needs to be calculated on such a segment, and then the solution of the equation can be derived recursively. Numerical examples demonstrate that it is a stable and efficient algorithm for the waveguide equations. This method can greatly reduces the requirement of internal memory and the amount of computation compared with the traditional algorithms.

**Keyword:** non-homogeneous waveguide equation, DtN map, piecewise uniform waveguide

## 1. Introduction

For two dimensional waveguide, in frequency domain, the governing equation is Helmholtz equation. The traditional algorithms for it lead to large linear systems in the calculation process, with time-consuming and memory-consuming. For homogeneous Helmholtz equation, some numerical methods based on DtN maps came out in 1996 and they have been improved in the practical waveguide calculation[1,2,5,6].

The waveguide propagation with inside light sources can be described by non-homogeneous Helmholtz equation. The numerical methods for the equation mainly are the traditional algorithms. Optical waveguides with piecewise uniform structure along the main propagation direction are widely used. Efficient algorithms for such waveguides are required in the design. In this paper, an efficient algorithm for non-homogeneous Helmholtz equation with piecewise uniform structure is derived, based on the modified DtN map of each uniform segment. The DtN map is computed numerically by a method using a difference approximation of 4[th] order to the second derivative for each segment. For segments with the same structure, the DtN map of each segment is the same, so only a DtN map needs to be calculated. With the computed DtN maps, the solution of the equation can be derived recursively. In section 2, non-homogeneous Helmholtz equation and the related boundary conditions are described. The new algorithm is introduced in section 3. Two numerical examples are given in section 4.

## 2. Non-homogeneous Helmholtz Equation

For transverse electric polarization with inside light source, the governing equation in frequency domain is

$$u_{zz} + u_{xx} + k_0^2 n^2(x,z) u = f(x,z) \qquad (1)$$

Where $u$ is the y-component of the electric field, $k_0$ is the free space wavenumber, $n(x,z)$ is the refractive index of the medium and $f(x,z)$ is the function of the inside light source.

We assume the waveguide axis is $z$ and the waveguide is $z$-invariant for $z < z_0$ and $z > z_m$ and piecewise $z$-invariant in $[z_0, z_m]$. Namely, we assume

---

[1] corresponding author, Email: huangjy@mail.buct.edu.cn

$$z_{-1} = -\infty, \quad z_0 < z_1 < \cdots < z_m, \quad z_{m+1} = \infty,$$

and

$$n(x,z) = n_j(x) \quad \text{for } z_{j-1} < z < z_j, \quad j = 0,1,\ldots,m+1$$

For an open structure, we use perfectly matched layer (PML) technique[3] to truncate unbounded domains. The waveguide problem with the inside light source and incident wave of the same time frequency becomes:

$$\begin{cases} \partial_z^2 u + Lu = f(x,z) \quad (x \in (-D,D), z \in (z_0, z_m)) \\ (L = \dfrac{1}{s(x)} \dfrac{\partial}{\partial x}(\dfrac{1}{s(x)} \dfrac{\partial}{\partial x}) + k_0^2 n^2(x,z)) \\ u(-D,z) = u(D,z) = 0 \\ (\partial_z u + i\sqrt{L(z_0^-)}\, u)|_{z=z_0} = 2i\sqrt{L(z_0^-)}\, u^+(x, z_0^-) \\ (\partial_z u - i\sqrt{L(z_m^+)}\, u)|_{z=z_m} = 0 \end{cases} \quad (2)$$

where $u^+$ is the incident wave, $u^-$ the reflected wave, and $s(x) = 1 + i\sigma(x)$, $\sigma(x) \geq 0$.

### 3. The new stable algorithm based on DtN maps

For the non-homogeneous Helmholtz equation, let the Riccati operator $Q$ and the fundamental solution operator $Y$, the functions $g(x,z)$ and $h(x,z)$ satisfy:

$$\begin{cases} u_z = Q(z)u + g \quad (z_0 < z < z_m) \\ Q(z_m) = i\sqrt{L(z_m^+)}, \quad g(z_m) = 0 \\ Y(z)u(x,z) + h = u(x, z_m) \quad (z_0 < z < z_m) \\ Y(z_m) = I, \quad h(z_m) = 0 \end{cases} \quad (3)$$

We set a DtN map $M_j$ and a vector valued function $s_j$ in the interval $(z_{j-1}, z_j)$ to satisfy

$$M_j \begin{bmatrix} u(x, z_{j-1}) \\ u(x, z_j) \end{bmatrix} + s_j = \begin{bmatrix} \partial_z u(x, z_{j-1}) \\ \partial_z u(x, z_j) \end{bmatrix} \quad (4)$$

$M_j$ is the same for each segment with the same waveguide structure. Denote $M_j = \begin{pmatrix} M_{11} & M_{12} \\ M_{21} & M_{22} \end{pmatrix}$ in block form, and $s_j = (s_1, s_2)^T$. Once $M_j$ and $s_j$ are computed, from (3) and (4), we can get a recurrence formula of $Q$, $Y$, $g$ and $h$:

$$\begin{cases} Q(z_{j-1}) = M_{11} + M_{12}(Q(z_j) - M_{22})^{-1} M_{21} \\ Y(z_{j-1}) = Y(z_j)(Q(z_j) - M_{22})^{-1} M_{21} \\ g(z_{j-1}) = s_1 + M_{12}(Q(z_j) - M_{22})^{-1}(s_2 - g(z_j)) \\ h(z_{j-1}) = h(z_j) + Y(z_j)(Q(z_j) - M_{22})^{-1}(s_2 - g(z_j)) \end{cases} \quad (5)$$

From $Q(z_m) = i\sqrt{L(z_m^+)}$, $Y(z_m) = I$, $g(z_m) = 0$, $h(z_m) = 0$, we can get $\{Q(z_k), Y(z_k), g(z_k), h(z_k)\}$ $(k = m-1,\cdots,0)$ recursively. If $Q(z_0), Y(z_0), g(z_0), h(z_0)$ have been obtained, from $\partial_z u(x, z_0) = Q(z_0)u(x, z_0) + g(z_0)$ and $(\partial_z u + i\sqrt{L(z_0^-)}\, u)|_{z=z_0} = 2i\sqrt{L(z_0^-)}\, u^+(x, z_0^-)$, we can get

$$u(x, z_0) = (Q(z_0) + i\sqrt{L(z_0^-)})^{-1}(2i\sqrt{L(z_0^-)}u^+(x, z_0^-) - g(z_0)) \tag{6}$$

Then, from $Y(z)u(x,z) + h = u(x, z_m)$ we can obtained $u(x, z_m)$ and $u(x, z_k)(k = 1, 2, \cdots, m-1)$.

Noting we only need compute the DtN map M of a segment for all segments with the same structure, the algorithm can significantly save computing time and memory.

We give an algorithm for computing the DtN map M in the following. We know there is the following approximation

$$\frac{1}{12}[f''(x_{j-1}) + 10 f''(x_j) + f''(x_{j+1})]$$
$$= \frac{1}{h^2}[f(x_{j-1}) - 2f(x_{j-1}) + f(x_{j-1})] + O(h^4) \tag{7}$$

For the interval $(z_{j-1}, z_j)$, we set

$$z_{j-1} = \xi_0 < \xi_1 < \cdots < \xi_q = z_j, \xi_k = \xi_0 + kh \tag{8}$$

Applying the approximation formula (7), we can get approximately:

$$\frac{1}{12}\begin{bmatrix} 10 & 1 & & \\ 1 & 10 & \ddots & \\ & \ddots & \ddots & 1 \\ & & 1 & 10 \end{bmatrix}_{q-1 \times q-1} \begin{bmatrix} u''(\xi_1) \\ u''(\xi_2) \\ \vdots \\ u''(\xi_{q-1}) \end{bmatrix}$$
$$= \frac{1}{h^2}\begin{bmatrix} -2 & 1 & & \\ 1 & -2 & \ddots & \\ & \ddots & \ddots & 1 \\ & & 1 & -2 \end{bmatrix}\begin{bmatrix} u(\xi_1) \\ u(\xi_2) \\ \vdots \\ u(\xi_{q-1}) \end{bmatrix} + \frac{1}{h^2}\begin{bmatrix} u(\xi_0) \\ 0 \\ \vdots \\ u(\xi_q) \end{bmatrix} - \frac{1}{12}\begin{bmatrix} u''(\xi_0) \\ 0 \\ \vdots \\ u''(\xi_q) \end{bmatrix} \tag{9}$$

Denoting

$$A = \frac{1}{12}\begin{bmatrix} 10 & 1 & & \\ 1 & 10 & \ddots & \\ & \ddots & \ddots & 1 \\ & & 1 & 10 \end{bmatrix}_{(q-1)\times(q-1)}, B = \begin{bmatrix} -2 & 1 & & \\ 1 & -2 & \ddots & \\ & \ddots & \ddots & 1 \\ & & 1 & -2 \end{bmatrix}$$

$$U = \begin{bmatrix} u(\xi_1) \\ u(\xi_2) \\ \vdots \\ u(\xi_{q-1}) \end{bmatrix}, \quad a_0 = \begin{bmatrix} 1 \\ 0 \\ \vdots \\ 0 \end{bmatrix}, \quad a_q = \begin{bmatrix} 0 \\ \vdots \\ 0 \\ 1 \end{bmatrix}$$

we have

$$AU_{zz} = \frac{1}{h^2}[a_0 u(x, z_{j-1}) + BU + a_q u(x, z_j)] - \frac{1}{12}[a_0 u''(x, z_{j-1}) + a_q u''(x, z_j)].$$

From (2) we have

$$U_{zz} + LU = F \qquad (10)$$

where $F = [f(x,\xi_1),\cdots,f(x,\xi_{q-1})]^T$. Using $u''(x,z_{j-1}) = -Lu(x,z_{j-1}) + f(x,z_{j-1})$ and $u''(x,z_j) = -Lu(x,z_j) + f(x,z_j)$, we can get

$$LU + \frac{1}{h^2}DU = -\frac{1}{h^2}[\hat{a}_0 u(x,z_{j-1}) + \hat{a}_q u(x,z_j)] - \frac{1}{12}[\hat{a}_0 Lu(x,z_{j-1}) + \hat{a}_q Lu(x,z_j)] + \frac{1}{12}[\hat{a}_0 f(x,z_{j-1}) + \hat{a}_q f(x,z_j)] + F \qquad (11)$$

where $D = A^{-1}B$, $\hat{a}_0 = A^{-1}a_0$, $\hat{a}_q = A^{-1}a_q$. Diagonalizing similarly the matrix $D$, we have

$$D = R\begin{bmatrix} \mu_1 & & & \\ & \mu_2 & & \\ & & \ddots & \\ & & & \mu_{q-1} \end{bmatrix} R^{-1} \qquad (12)$$

where $\mu_k = 12\dfrac{-1+\cos(\frac{k\pi}{q})}{5+\cos(\frac{k\pi}{q})}$, $k = 1,\ldots,q-1$, and $R = \left(\sin(\dfrac{jk\pi}{q})\right)$.

Finally, we get $q-1$ independent ordinary differential equations:

$$Lw_k + \frac{1}{h^2}\mu_k w_k = -\frac{1}{h^2}[\alpha_k u(x,z_{j-1}) + \beta_k u(x,z_j)] - \frac{1}{12}[\alpha_k Lu(x,z_{j-1}) + \beta_k Lu(x,z_j)] + \frac{1}{12}[\alpha_k f(x,z_{j-1}) + \beta_k f(x,z_j)] + \zeta_k \quad (k = 1,\cdots,q-1) \qquad (13)$$

where

$$\begin{bmatrix} w_1 \\ w_2 \\ \vdots \\ w_{q-1} \end{bmatrix} = R^{-1}U, \quad \begin{bmatrix} \alpha_1 \\ \alpha_2 \\ \vdots \\ \alpha_{q-1} \end{bmatrix} = R^{-1}a_0, \quad \begin{bmatrix} \beta_1 \\ \beta_2 \\ \vdots \\ \beta_{q-1} \end{bmatrix} = R^{-1}a_q, \quad \begin{bmatrix} \zeta_1 \\ \zeta_2 \\ \vdots \\ \zeta_{q-1} \end{bmatrix} = R^{-1}F, \quad F = \begin{bmatrix} f(x,\xi_1) \\ f(x,\xi_2) \\ \vdots \\ f(x,\xi_{q-1}) \end{bmatrix} \qquad (14)$$

Solving (11), we can get $w_k(x)$ and then $U$.

If we can compute $\partial_z u(x,z_{j-1})$ and $\partial_z u(x,z_j)$ from any given $u(x,z_{j-1})$ and $u(x,z_j)$, then from (4), we can get the map $M_j$ and the vector $s_j$.

Using $u''(x,z) = -Lu(x,z) + f(x,z)$, we have the following approximations with error $O(h^4)$

$$\partial_z u(x, z_{j-1}) \approx \frac{u(x, \xi_2) - u(x, z_{j-1})}{2h}$$
$$- h \frac{[-Lu(x, z_{j-1}) + f(x, z_{j-1})] + 2[-Lu(x, \xi_1) + f(x, \xi_1)]}{3} \quad (15)$$

$$\partial_z u(x, z_j) \approx \frac{u(x, z_j) - u(x, \xi_{q-2})}{2h}$$
$$+ h \frac{2[-Lu(x, \xi_{q-1}) + f(x, \xi_{q-1})] + [-Lu(x, z_j) + f(x, z_j)]}{3} \quad (16)$$

If we discrete $[-D, D]$ with $-D = x_0 < x_1 < \cdots < x_{N+1} = D$, $x_n = x_0 + nh_x$, ($n = 0, 1, \cdots, N+1$), from (4), we get approximately:

$$\begin{bmatrix} \partial_z u(x_1, z_{j-1}) \\ \vdots \\ \partial_z u(x_N, z_{j-1}) \\ \partial_z u(x_1, z_j) \\ \vdots \\ \partial_z u(x_N, z_j) \end{bmatrix} = \begin{bmatrix} M_{11} & M_{12} \\ M_{21} & M_{22} \end{bmatrix} \begin{bmatrix} u(x_1, z_{j-1}) \\ \vdots \\ u(x_N, z_{j-1}) \\ u(x_1, z_j) \\ \vdots \\ u(x_N, z_j) \end{bmatrix} + \begin{bmatrix} s_1(x_1) \\ \vdots \\ s_1(x_N) \\ s_2(x_1) \\ \vdots \\ s_2(x_N) \end{bmatrix} \quad (17)$$

Here $M_{ij}(i, j = 1, 2)$ are $N \times N$ matrices. The equations (13) becomes:

$$LW_k + \frac{1}{h^2} \mu_k W_k = -\frac{1}{h^2}(\alpha_k U_{j-1} + \beta_k U_j) - \frac{1}{12}[\alpha_k L U_{j-1} + \beta_k L U_j] + \frac{1}{12}[\alpha_k F_{j-1} + \beta_k F_j] + Z_k \quad (18)$$

where

$$W_k = \begin{pmatrix} w_k(x_1) \\ \vdots \\ w_k(x_N) \end{pmatrix}, U_{j-1} = \begin{pmatrix} u(x_1, z_{j-1}) \\ \vdots \\ u(x_N, z_{j-1}) \end{pmatrix}, U_j = \begin{pmatrix} u(x_1, z_j) \\ \vdots \\ u(x_N, z_j) \end{pmatrix},$$

$$F_{j-1} = \begin{pmatrix} f(x_1, z_{j-1}) \\ \vdots \\ f(x_N, z_{j-1}) \end{pmatrix}, F_j = \begin{pmatrix} f(x_1, z_j) \\ \vdots \\ f(x_N, z_j) \end{pmatrix}, Z_k = \begin{pmatrix} \zeta_k(x_1) \\ \vdots \\ \zeta_k(x_N) \end{pmatrix}$$

and $L$ is a $N \times N$ matrix which approximates the operator $L$ in (2). First, we put $\begin{bmatrix} u(x_1, z_{j-1}) \\ \vdots \\ u(x_N, z_{j-1}) \end{bmatrix} = \begin{bmatrix} 0 \\ \vdots \\ 0 \end{bmatrix}$ and

$\begin{bmatrix} u(x_1, z_j) \\ \vdots \\ u(x_N, z_j) \end{bmatrix} = \begin{bmatrix} 0 \\ \vdots \\ 0 \end{bmatrix}$. For every $k$, we can get $W_k$ from (18) and then $U(x_j) = R[w_1(x_j), \cdots, w_{q-1}(x_j)]^T$ from (14). Next, we can compute $\begin{bmatrix} \partial_z u(x_1, z_{j-1}) \\ \vdots \\ \partial_z u(x_N, z_{j-1}) \end{bmatrix}$ and $\begin{bmatrix} \partial_z u(x_1, z_j) \\ \vdots \\ \partial_z u(x_N, z_j) \end{bmatrix}$ from (15) and (16) and then give

$$\begin{bmatrix} s_1(x_1,z_{j-1}) \\ \vdots \\ s_1(x_N,z_{j-1}) \end{bmatrix} \text{ and } \begin{bmatrix} s_2(x_1,z_{j-1}) \\ \vdots \\ s_2(x_N,z_{j-1}) \end{bmatrix} \text{ from (17)}.$$

After that, similarly, if we take $[u(x_1,z_{j-1}),\cdots,u(x_N,z_{j-1}),u(x_1,z_j),\cdots,u(x_N,z_j)]^T$ as each column of the unit matrix of $(2N)\times(2N)$ respectively, we can get each column of the matrix $M_j$.

## 4. Example

We use the developed algorithm to compute an example of COST 268 modeling task[4] with an addition of non-homogeneous term $f(x,z)$ but without incident wave.

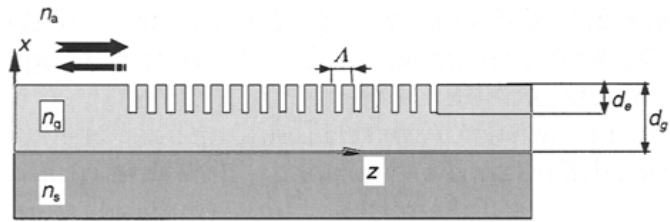

**Fig4-1** Bragg grating in a planar waveguide

It is about a high-contrast optical waveguide with a deeply etched short Bragg grating. The Bragg grating is composed of 20 rectangular grooves. The widths of the "tooth" and "groove" both are equal to 0.215μm. The numerical results are shown in the following figure.

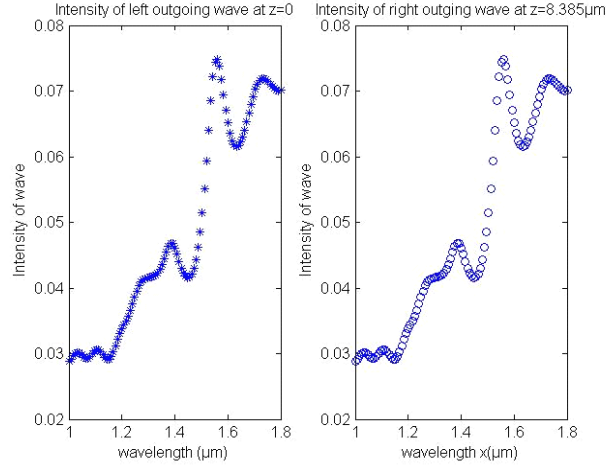

**Fig4-2** The intensity of left and right outgoing waves for waveguide with inside light source $f(x,z)=\cos(\dfrac{\pi x}{2D})\sin(\dfrac{\pi z}{0.215})$ at the center tooth and wavelengths from 1 μm to 1.8 μm

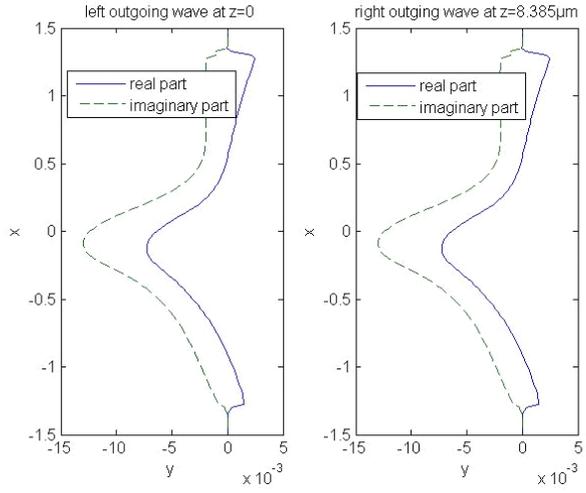

**Fig4-3** Left and right outgoing waves for waveguide with inside light source $f(x,z) = \cos(\frac{\pi x}{2D})\sin(\frac{\pi z}{0.215})$ at the center tooth and wavelength 1.8 μm

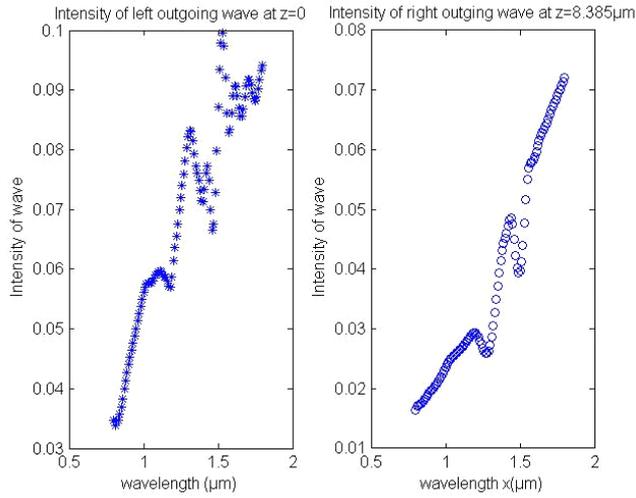

**Fig4-4** The intensity of left and right outgoing waves for waveguide with inside light source $f(x,z) = \cos(\frac{\pi x}{2D})\sin(\frac{\pi z}{0.215})$ at the 2$^{nd}$ tooth from left and wavelengths from 0.8 μm to 1.8 μm

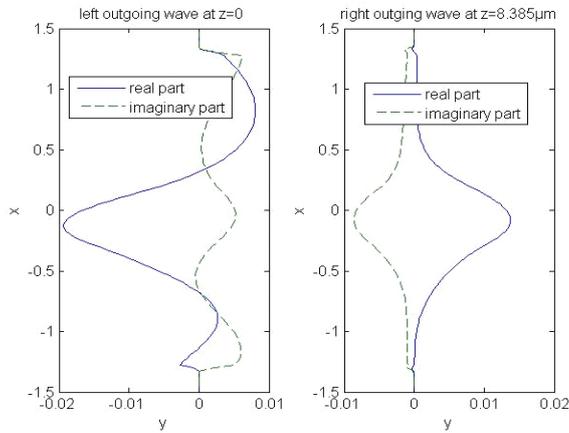

**Fig4-5** Left and right outgoing waves for waveguide with inside light source

$$f(x,z) = \cos(\frac{\pi x}{2D})\sin(\frac{\pi z}{0.215})$$ at the 2nd tooth from left and wavelength 1.8 μm

## 5. Conclusions

Using the DtN map, for non-homogeneous Helmholtz equation with a piecewise uniform structure, the new algorithm can greatly reduce the computing time and memory, because, for segments with the same piecewise uniform structure, the DtN map of each segment is the same, the calculation of the DtN map can be reduced to solving some independent ordinary differential equations, and the solution of the non-homogeneous Helmholtz equation can be computed by marching on in the main propagation direction. Similarly, if solving the equation (11) directly, the algorithm is also efficient for waveguides with general periodic structure or partial periodic structure because the DtN map of each periodic section is the same too. The numerical examples demonstrate that the algorithm is a stable and efficient.